\documentclass[11pt,reqno]{amsart}
\usepackage{amscd,graphicx,psfrag}

\newcommand{\p}{\partial}

\newcommand{\Z}{\mathbb Z}

\renewcommand{\L}{\mathcal L}

\newcommand{\M}{\mathcal M}
\newcommand{\R}{\mathcal R}
\newcommand{\E}{\mathcal E}

\renewcommand{\phi}{\varphi}

\newcommand{\tr}{\operatorname{tr}\,}

\newcommand{\ad}{\operatorname{ad}}

\newcommand{\Fix}{\operatorname{Fix}\,}
\newcommand{\Ad}{\operatorname{Ad}}

\newcommand{\End}{\operatorname{End}}

\newcommand{\N}{\mathcal N}
\newcommand{\si}{\sigma}
\newcommand{\ov}{\overline}

\newtheorem{theorem}{Theorem}[section]
\newtheorem{lemma}[theorem]{Lemma}
\newtheorem{proposition}[theorem]{Proposition}
\newtheorem{corollary}[theorem]{Corollary}

\theoremstyle{definition}
\newtheorem{remark}[theorem]{Remark}

\title{On real moduli spaces over M--curves}
\author{Nikolai Saveliev}
\address{Department of Mathematics\newline\indent
University of Miami \newline\indent PO Box 249085
\newline\indent Coral Gables, FL 33124}
\email{\rm{saveliev@math.miami.edu}}
\author{Shuguang Wang}
\address{Deapartment of Mathematics\newline\indent
University of Missouri \newline\indent
202 Mathematical Sciences Bldg \newline\indent
Columbia, MO 65211}
\email{\rm{sw@math.missouri.edu}}
\thanks{The first author was partially supported by NSF Grant
DMS 0305946. The second author was patially supported by the 
University of Missouri Research Board Grant.}

\begin{document}

\begin{abstract}
Let $F$ be a genus $g$ curve and $\sigma: F \to F$ a real structure
with the maximal possible number of fixed circles. We study the 
real moduli space $\N' = \Fix (\sigma^{\#})$ where $\sigma^{\#}: 
\N \to \N$ is the induced real structure on the moduli space $\N$
of stable holomorphic bundles of rank 2 over $F$ with fixed 
non-trivial determinant. In particular, we calculate $H^* 
(\N',\mathbb Z)$ in the case of $g = 2$, generalizing Thaddeus' 
approach to computing $H^* (\N,\mathbb Z)$. 
\end{abstract}

\maketitle


\section{Introduction} 

In their influential paper \cite{AB}, Atiyah and Bott used two--dimensional 
Yang--Mills theory to compute cohomology of the moduli space $\N$ of stable 
holomorhic bundles over a Riemann surface $F$.  In essence, the computation 
was inspired by the idea that the Yang--Mills functional in this dimension 
should be an equivariantly perfect Morse--Bott function (this was proved 
later in full generality by Daskalopoulos \cite{Da}). 

A real structure $\sigma: F\to F$ induces a real structure $\sigma^{\#}$ on 
the moduli space $\N$. Understanding the structure of the real moduli space
$\N' = \Fix(\sigma^{\#})$ is an important but subtle problem. It has been 
addressed most recently in a series of papers by N.-K.~Ho and C.-C.~M.~Li,
the latest being \cite{HL}. They mainly treat the case when the real 
structure $\sigma$ has no fixed points, which leads them to the study of 
the Yang--Mills functional on the non--orientable surface $F/\sigma$.

In this paper we consider the other extreme,  when $\sigma$ has the maximal
possible number of fixed circles.  In this case,  the pair  $(F,\sigma)$ is 
usually referred to as an $M$--curve.  More specifically, we work with the 
moduli space $\N$ of stable holomorphic bundles of rank 2 over $F$ with
fixed non--trivial determinant, and with the associated real moduli space 
$\N'$. Instead of the infinite dimensional Yang--Mills functional we utilize 
a finite dimensional Morse--Bott function as in Thaddeus \cite{T}. The 
perfection of the function in Thaddeus' paper was suggested by the work of 
Frankel \cite{F} while the perfection of ours is suggested by Duistermaat's 
paper \cite{D}.

For technical reasons, in this paper we will only take up the case of genus 
two $M$--curves, hoping to address the case of higher genera elsewhere. For
$M$--curves of genus two, the second author \cite{W1} described $\N'$ 
algebraically as the intersection of two quadrics in the five--dimensional 
real projective space; however, it proved to be rather difficult to extract 
any further information about the topology of $\N'$ from that description. 
The Morse theoretic approach of this paper, on the other hand, gives a 
complete calculation of the integral cohomology of $\N'$.

\begin{theorem}\label{T:main}
Let $(F,\sigma)$ be a genus two $M$--curve and $\N'$ the real moduli space 
of stable holomorphic bundles of rank 2 over $F$ with fixed non--trivial
determinant. Then, at the level of graded abelian groups,  there is an 
isomorphism $H^* (\N',\Z) = H^* (S^1 \times S^1 \times S^1,\Z)$.
\end{theorem}

It should be mentioned that the above isomorphism does not hold at the 
level of cohomology rings, see Remark \ref{R:rings}.

Parts of this project were accomplished while we attended the Summer 2006 
Session of the Park City/IAS Mathematics Institute. We express our 
appreciation to the organizers for providing a stimulating environment.


\section{Real moduli spaces over curves}
Let $F$ be a compact surface of genus $g \ge 2$ and fix a point $p \in F$. 
Denote by $\N$ the moduli space of stable holomorphic bundles $\E \to F$ 
of rank 2 with determinant $\L_p^{-1}$, where $p \in F$ is viewed as a 
divisor. Then $\N$ is a smooth complex manifold of real dimension $6g - 6$, 
modeled on the complex vector space $H^1 (F,\End \E)$ by the deformation 
theory.

Let $\sigma: F \to F$ be a real structure on $F$ whose real part $F' = 
\Fix(\sigma)$ is non--empty. Then $F'$ contains at least one circle, 
but may contain as many as $g + 1$ circles. In the latter case, the pair 
$(F,\sigma)$ is called an $M$--curve. Choose $p \in F'$ then $\sigma: F 
\to F$ induces an involution ${\sigma^{\#}}: \N \to \N$ by the formula 
\[
\sigma^{\#}\,[\E] = [\,\sigma^*\,\ov{\E}\,], 
\]
where $\ov{\E}$ stands for the complex conjugate of $\ov{\E}$. Note that 
since $\sigma: F \to F$ is orientation reversing, $\sigma^*{\L_p} = 
\ov{\L}_p$, thus making the complex conjugation in the above formula 
necessary. 

Note that the linear map $H^1 (F,\End \E) \to H^1(F,\End (\si^*\ov{\E}))$ 
induced by $\sigma$ is a complex conjugation, therefore, the involution 
$\sigma^{\#}: \N \to \N$ is anti-holomorhic and hence is a real structure 
on $\N$. It follows that the real moduli space $\N' = \Fix (\sigma^{\#})$ 
is a closed smooth manifold of dimension $3g - 3$. 

\begin{proposition}
The manifold $\N'$ is orientable. 
\end{proposition}

\begin{proof}
According to \cite{AB}, the moduli space $\N$ is simply connected and spin. 
Therefore, the real structure $\sigma^{\#}$ must be compatible with the 
unique spin structure on $\N$ in the sense of \cite{W2}. The main result 
in \cite{W2} then applies to infer that $\N'$ is orientable. Note that 
$\sigma^{\#}$ need not preserve the spin structure in the usual sense, 
since it is orientation reversing when the complex dimension of $\N$ is 
odd, that is, when $g$ is even.
\end{proof}


\section{Representation varieties}
Let $F_0$ be the surface $F$ punctured at $p \in F$. The theorem of 
Narasimhan and Seshadri \cite{NS} can be used to identify $\N$ with the 
representation variety $\M$ which consists of the conjugacy classes of 
$SU(2)$ representations of $\pi_1 (F_0) = \pi_1 (F_0, x_0)$ sending 
$[\,\p F_0\,] \in \pi_1 (F_0)$ to $-1 \in SU(2)$  (the latter condition 
does not dependent on the choice of basepoint because $-1$ is a central 
element in $SU(2)$).

Given a real structure $\sigma: F\to F$ with non--empty $F'= \Fix(\sigma)$ 
and $p \in F'$ as above, choose a basepoint $x_0 \in F' \cap F_0$. Then we 
have an induced involution $\sigma_*: \pi_1 (F_0) \to \pi_1 (F_0)$, which 
in turn induces an involution $\sigma^*: \M \to \M$ by the formula
\[
\sigma^*\,[\alpha] = [\,\alpha\circ\sigma_*\,].
\]
That $\sigma^*$ is well defined follows from the fact that $\sigma_*\,
[\,\p F_0\,] = [\,\p F_0\,]^{-1} = -1$. 

\begin{lemma}
Let $\phi: \M \to \N$ be the Narasimhan--Seshadri diffeomorphism then 
the following diagram commutes
\[
\begin{CD}
\M @> \sigma^* >> \M \\
@V \phi VV @V \phi VV \\
\N @> \sigma^{\#} >> \N
\end{CD}
\]
\end{lemma}

\medskip

\begin{proof}
The Narasimhan--Seshadri correspondence assigns to every $[\alpha] \in 
\M$ a stable holomorphis bundle $\E_{\alpha} \to F$ as follows. Let 
$\tilde F_0 \to F_0$ be the universal covering space of $F_0$ and 
consider the holomorphic bundle $\E \to F_0$ with
\[
\E\; =\; \tilde F_0\, \times_{\pi_1(F_0)}\, \mathbb C^2, 
\]
where $\pi_1 (F_0)$ acts on $\mathbb C^2$ via $\alpha: \pi_1 (F_0) \to 
SU(2)$. Obviously, $\det \E$ is trivial. Since the holonomy of $\alpha$ 
along the loop $\p F_0$ is fixed, we can trivialize $\E$ near the boundary 
$\p F_0$. Glue $D^2 \times \mathbb C^2$ in using the transition function 
$z^{-1}$ along $\p F_0$. The result is the stable bundle $\E_{\alpha} 
\to F$ with the determinant $\det\E_{\alpha} = \L^{-1}_p$ yielding the 
Narasimhan--Seshadri correspondence. 

For any given $[\alpha] \in \M$ we have $\sigma^{\#}\,[\E_{\alpha}] = 
[\,\sigma^*\,{\ov \E_{\alpha}}\,] = [\,\E_{\sigma^* \ov \alpha}\,]$, 
where $\ov \alpha$ is the complex conjugate of $\alpha$. However, for 
any matrix 
\[
A\; =\;
\left[\begin{array}{cc}
a & b\\
- \ov b & \ov a
\end{array}
\right]
\;\in\; SU(2),
\]

\noindent
its complex conjugate 

\[
\ov A \; =\;
\left[\begin{array}{cc}
\ov a & \ov b\\
- b & a
\end{array}
\right]
\;=\;
\left[\begin{array}{cc}
0  & 1\\
-1 & 0
\end{array}
\right]
\;A\; 
\left[\begin{array}{cc}
0  & 1\\
-1 & 0
\end{array}
\right]^{-1}
\]

\medskip\noindent
is the matrix conjugate of $A$ via 
\[
j\; =\; 
\left[\begin{array}{cc}
0  & 1\\
-1 & 0
\end{array}
\right]
\;\in\; SU(2)
\] 

\noindent
(we use the standard identification between $SU(2)$ matrices and unit 
quaternions). This means that $\ov \alpha$ and $\alpha$ are conjugate 
representations, and therefore the bundles $\E_{\sigma^* \ov \alpha}$ 
and $\E_{\sigma^* \alpha}$ are isomorphic.
\end{proof}

Denote by $\M'$ the fixed point set of the involution $\sigma^*: \M \to 
\M$ then the above lemma implies that $\M'$ and $\N'$ are diffeomorphic.

\begin{corollary}\label{C:orient}
The moduli space $\M'$ is a smooth closed orientable manifold of 
dimension $3g - 3$.
\end{corollary}

In conclusion, note that $\M$ is a symplectic manifold with symplectic 
form $\omega: H^1 (F;\ad\rho)\,\otimes\,H^1 (F;\ad\rho)\to \mathbb R$ 
given by $\omega(u,v) = - 1/2\,\tr (u\,\cup\,v)\,[F]$, see Goldman 
\cite{G}.

\begin{lemma}\label{L:anti}
The map $\sigma^*: \M \to \M$ is a real sructure with respect to the 
symplectic form $\omega$, that is, $\sigma^* \omega = - \omega$. 
\end{lemma}

\begin{proof}
This result is obtained by the following straightforward calculation\,:
\begin{multline}\notag
(\sigma^* \omega) (u,v) =
\omega (\sigma^* u,\sigma^* v) = -1/2\,\tr (\sigma^* u\,\cup\,\sigma^* v)
\,[F] \\ = -1/2\,\tr (\sigma^* (u\,\cup\,v))\,[F] = -1/2\,\tr (u\,\cup\,v)
\,\sigma_*[F] = -\omega(u,v).
\end{multline} 
In the last equality, we used the fact that the map $\sigma: F \to F$ is 
orientation reversing. 
\end{proof}


\section{The case of $g = 2$}
\begin{figure}[!ht]
\psfrag{p}{$p$}
\psfrag{c}{$c$}
\psfrag{a1}{$a_1$}
\psfrag{a2}{$a_2$}
\psfrag{b1}{$b_1$}
\psfrag{b2}{$b_2$}
\centering
\includegraphics{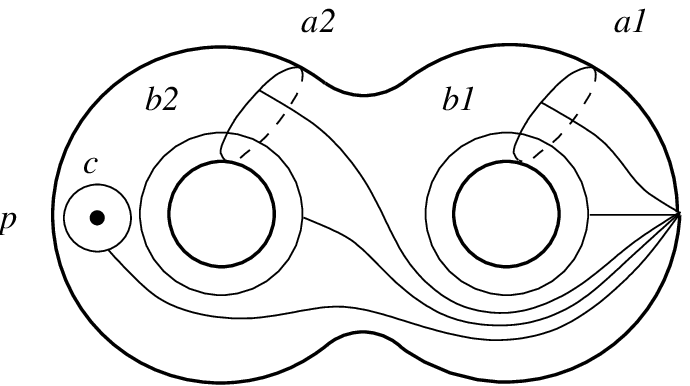}
\caption{}
\label{fig2}
\end{figure}

Let $(F,\sigma)$ be a genus two $M$--curve embedded in $\mathbb R^3$ as 
shown in Figure \ref{fig2}, with $\sigma: F \to F$ acting as the reflection 
in the horizontal plane fixing the obvious three circles.  Let $F_0$ be 
a once punctured surface $F$, the puncture $p \in F$ being a real point
on the left circle of $F'$  (the other two options for positioning $p$ 
will be treated in Section \ref{S:proof}). Then $\pi_1 (F_0)$ is the free 
group
\[
\pi_1 (F_0) = \langle\,a_1,b_1,a_2,b_2\,|\qquad\rangle
\]
whose generators $a_1$, $a_2$, $b_1$, $b_2$ are shown in the picture. 
Also shown is the curve $c = [a_1,b_1] \cdot [a_2,b_2]$. All of the 
above curves are oriented clockwise with respect to the projection 
shown in Figure \ref{fig2}. 

\bigskip

The moduli space $\M$ can now be described as follows. Let us consider 
the smooth map $\mu: SU(2)^4 \to SU(2)$ given by $\mu (A_1,B_1,A_2,B_2) 
= - [A_1,B_1]\cdot [A_2,B_2]$. Since $1\in SU(2)$ is a regular value of 
$\mu$, its preimage $\mu^{-1} (1)$ is a smooth manifold. It contains no 
reducibles, that is, 4--tuples $(A_1,B_1,A_2,B_2)$ whose entries commute 
with each other, because the latter would contradict the equation 
$[A_1,B_1] \cdot [A_2,B_2] = -1$.  Therefore, the action of $SO(3) 
= SU(2)/\pm 1$ on $\mu^{-1}(1)$ by conjugation is free, and its quotient 
space is then the smooth six--dimensional manifold $\M$.
A straightforward calculation shows that the induced action 
$\sigma: \pi_1 (F_0) \to \pi_1 (F_0)$ is given on the generators by the 
formulas
\begin{alignat*}{1}
& \sigma_* (a_1) = a_1,\quad\sigma_*(b_1) = b_1^{-1}, \\ 
&\sigma_* (a_2) = b_1^{-1} c\,b_2\,a_2\,b_2^{-1} b_1,\quad
\sigma_* (b_2) = b_1^{-1} b_2^{-1} b_1.
\end{alignat*}

\noindent
In practical terms, $\M$ consists of the conjugacy classes  $[A_1,B_1,A_2,
B_2]$ of quadruples $(A_1,B_1,A_2,B_2)$ such that $[A_1,B_1] \cdot [A_2,
B_2] = -1$, and the real structure $\sigma^*: \M \to \M$ is given by
\[
\sigma^* [A_1,B_1,A_2,B_2] = [B_1\,A_1\,B_1^{-1},\;B_1^{-1},\;
-B_2\,A_2\,B_2^{-1},\;B_2^{-1}].
\]

\noindent
As we already know, the fixed point set $\M'$ of $\sigma^*: \M \to \M$ is 
a smooth orientable manifold of dimension 3. 


\section{The function}
Let $f: \M \to \mathbb R$ be the function on the moduli space $\M$ defined 
by the formula 
\[
f\,([A_1,B_1,A_2,B_2]) = \tr (B_1)/2.
\]
This is a smooth function whose range is $[-1,1]$, and $f^{-1}(-1,1)$ is
acted upon by $S^1$ as follows.  After conjugation if necessary,  every 
element of $f^{-1}(-1,1)$ can be written in the form $[A_1,B_1,A_2,B_2]$ 
with $B_1 = e^{i\beta}$ and $0 < \beta < \pi$. The circle action is then 
given by $e^{i\phi}: [A_1,B_1,A_2,B_2] \mapsto [A_1\,e^{i\phi},B_1,A_2,
B_2]$.

\begin{lemma}
This is a well defined action on $f^{-1} (-1,1)$.
\end{lemma}

\begin{proof}
Any other choice of representative in $[A_1,B_1,A_2,B_2]$ with $B_1 = 
e^{i\gamma}$ and $0 < \gamma < \pi$, will have the property that 
$e^{i\gamma} = g e^{i\beta} g^{-1}$ for some $g \in SU(2)$. Therefore, 
$\gamma = \beta$ and $g$ is a unit complex number. Because of that, 
$(g A_1 g^{-1}) e^{i\phi} = g (A_1 e^{i\phi}) g^{-1}$, making the 
action well defined. Of course, if $B_1 = \pm 1$, the conjugating 
element $g$ can be any $SU(2)$ matrix, and the above argument fails.
\end{proof}

According to Goldman \cite{G}, the above circle action preserves the 
symplectic form $\omega$ on $\M$ and $\arccos f$ is its moment map 
(up to a factor of $i$).

It is clear that $f\circ\sigma^* = f$ hence $f: \M \to \mathbb R$ is 
invariant with respect to $\sigma^*: \M \to \M$. The restriction of $f$ 
to $\M'$ will be denoted by $f': \M'\to \mathbb R$. The range of $f'$ 
is again $[-1,1]$. The circle action on $f^{-1}(-1,1)$ is not defined 
on $(f')^{-1}(-1,1)$. Nevertheless, it gives rise to the residual 
$\Z/2$ action $r: \M' \to \M'$ defined on the entire moduli space $\M'$ 
(and in fact on the full moduli space $\M$) by the formula
\begin{equation}\label{E:res}
r\,([A_1,B_1,A_2,B_2]) = [-A_1,B_1,A_2,B_2]. 
\end{equation}


\section{The critical submanifolds of $f'$}
Thaddeus \cite{T} proved that $f: \M \to \mathbb R$ is a Morse--Bott 
function and described its critical submanifolds. These critical 
submanifolds are acted upon by $\sigma^*: \M \to \M$, and the fixed point 
sets of this action are exactly the critical submanifolds of $f': 
\M' \to \mathbb R$. These are of two types.

The first type is comprised of the manifolds $S_1'$, $S_3' \subset \M'$ on 
which $f'$ achieves its absolute minimum and maximum. To be 
precise, the absolute minimum $S_1 = f^{-1}(-1)$ of $f: \M \to \mathbb R$ 
is a copy of $SU(2)$ consisting of quadruples $[A_1,-1,i,j]$ parametrized 
by $A_1 \in SU(2)$. The action 
\[
\sigma^*[A_1,-1,\,i,\,j] = [A_1,-1,\,i,-j] = [i\,A_1\,i^{-1},-1,\,i,\,j] 
\]
corresponds to the map $\Ad_{\,i}: SU(2) \to SU(2)$ sending $A_1$ to 
$i\,A_1\,i^{-1}$. The fixed point set of this action, which is $S_1'$, is 
the circle consisting of quadruples $[e^{i\phi},-1,\,i,\,j]$ parametrized 
by $e^{i\phi}$. Similarly, $S_3 = f^{-1}(1)$, where $f: \M \to \mathbb R$ 
achieves its absolute maximum, gives rise to the circle $S_3'$ 
consisting of quadruples $[e^{i\psi},1,\,i,\,j]$.

\begin{lemma}
The critical submanifolds $S_1'$ and $S_3'$ of $f'$ are non-degenerate in 
the Morse--Bott sense, and their indices are respectively 0 and 2. 
\end{lemma}

\begin{proof}
The Hessian of $f$ is negative definite on the normal bundle of $S_3
\subset \M$. But then its restriction to the normal bundle of $S_3' 
\subset \M'$, which is the Hessian of $f'$, is also negative definite. 
In particular, $S_3'$ is non-degenerate, and its Morse--Bott index equals 
the codimension of $S_3'$ in $\M'$, which is 2. The argument for $S_1'$ 
is analogous.
\end{proof}

Within $f^{-1}(-1,1)$, the critical points of $f: \M \to \mathbb R$ 
coincide with those of the moment map $\arccos f$, hence with the fixed 
points of the circle action on $f^{-1}(-1,1)$. They form the submanifold 
$S_2 \subset \M$ which is a 2--torus consisting of quadruples $[j,\,i,\,
z,\,w]$ parametrized by $z$, $w \in \mathbb C$ such that $|z| = |w| = 1$. 
Note that $f = 0$ on $S_2$. One can easily see that the fixed point set 
$S_2'$ of the involution $\sigma^*: S_2 \to S_2$ given by
\[
\sigma^* [j,\,i,\,z,\,w] = [-j,-i,-z,\,\bar w] = [j,\,i,-\bar z,\,w]
\]
consists of the two circles $[j,\,i,\,\pm i,\,w]$ with $|w| = 1$.

\begin{lemma}
The critical submanifold $S_2'$ of $f'$ is non-degenerate in the 
Morse--Bott sense, and its index is equal to 1. 
\end{lemma}

\begin{proof}
According to Frankel \cite{F}, the moment map $\arccos f$ is a Morse--Bott 
function on $f^{-1}(-1,1)$, and hence so is the function $f$. The involution 
$\M \to \M$ given by $[A_1,B_1,A_2,B_2] \mapsto [A_1,-B_1,A_2,B_2]$ changes 
sign of $f$. Since $S_2$ is connected, the index of $f$ is half the rank of 
the normal bundle of $S_2 \subset \M$, or 2. The involution $\sigma^*$ is 
antisymplectic, see Lemma \ref{L:anti}, hence we can apply 
\cite[Proposition 2.2]{D} to the moment map $\arccos f$. It tells us that 
$\arccos f'$ and hence $f'$ are Morse--Bott on $S_2'$, and the index of 
$f'$ is half that of $f$.
\end{proof}

\begin{lemma}\label{L:action}
The residual involution $r: \M' \to \M'$ acts as the $180^{\circ}$ rotation 
on each of the circles $S_1'$ and $S_3'$, and acts trivially on $S_2'$. 
\end{lemma}

\begin{proof}
The circle $S_3'$ consists of the quadruples $[e^{i\psi},1,i,j]$ acted upon
by $r$ as $[e^{i\psi},1,i,j] \to [-e^{i\psi},1,i,j]$. The case of $S_1'$ is 
completely similar. The manifold $S_2'$ consists of the quadruples $[j,i,
\pm i,w]$ where $w$ is a unit complex number. The action of $r$ is given by 
$[j,i,\pm i,w]\to [-j,i,\pm i,w] = [j,i,\pm i,w]$, after conjugating by $i$.
\end{proof}


\section{The Morse--Bott spectral sequence}
\begin{figure}[!ht]
\psfrag{Z}{$\Z$}
\psfrag{Z2}{$\Z^2$}
\psfrag{0}{$0$}
\psfrag{p}{$p$}
\psfrag{q}{$q$}
\centering
\includegraphics{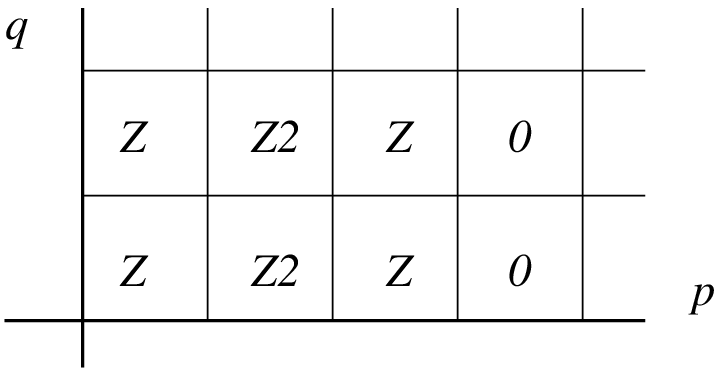}
\caption{}
\label{fig3}
\end{figure}

As we have seen, the critical set of the Morse--Bott function $f': \M'\to 
\mathbb R$ is a disjoint union $S_1'\,\cup\, S_2'\,\cup\, S_3'$, where 
the index of $S_p'$ is equal to $p - 1$, and the restriction of $f'$ 
to each of the $S_p'$ is constant, $f' (S_p') = p - 2$ for 
$p = 1, 2, 3$. Let $x_j = j - 3/2$ for $j = 0, 1, 2, 3$, and consider 
the filtration 
\[
\emptyset\; =\; U'_0\quad \subset\quad U'_1\quad \subset\quad 
U'_2\quad \subset\quad U'_3\; =\; \M'
\]
of $\M'$ by the open sets $U_j' = (f')^{-1}
(x_0,x_j)$. The Morse Lemma implies that, up to homotopy, $U'_j$ is a 
complex obtained from $U'_{j-1}$ by attaching, along its boundary, the 
disc bundle associated to the negative normal bundle over 
$S'_j$ whose fibers are the negative definite subspaces of the 
Hessian of $f'$ on $S'_j$. The $E_1$ term of the Morse--Bott spectral 
sequence associated with this filtration is of the form shown in Figure 
\ref{fig3}, where the $(p-1)$st column represents $H^* (S_p')$ for
$p = 1, 2, 3$. This spectral sequence converges to $H^*(\M')$.

\bigskip

Observe that the filtration is preserved by the residual $\Z/2$ action 
$r: \M' \to \M'$ hence $r$ induces an automorphism $r^*$ of the above 
spectral sequence. Our next task will be to compute the differentials, 
of which only $d_1$ and $d_2$ are not obviously zero. 


\section{Differentials $d_1$}
The differential $d_1: E_1^{0,0} \to E_1^{1,0}$ must vanish because $H^0 
(\M') \ne 0$. In particular, we conclude that $H^0(\M') = \Z$ and hence 
$\M'$ is connected. Similarly, the differential $d_1: E_1^{1,1} \to 
E_1^{2,1}$ must vanish because $\M'$ is orientable, see Corollary 
\ref{C:orient}, so $H^3 (\M') = \Z$. In fact, we will show that the 
remaining differentials $d_1$ also vanish but this will take some effort. 

Let us first consider the differential $d_1: E_1^{1,0}\to E_1^{2,0}$. It 
is represented as the composition 
\begin{equation}\label{E:thom1}
H^0 (S_2')\; \to\; H^1 (U_2',U_1')\; \to\; H^1 (U_2')\; \to\; 
H^2 (U_3',U_2')\; \to\; H^0 (S_3').
\end{equation}
Here, the second and third arrows are portions of the long exact sequences 
of the respective pairs hence they commute with the automorphism $r^*$. 
The first and the last arrows are Thom isomorphisms, therefore, their 
behavior with respect to $r^*$ is determined by how $r$ acts on $S_2'$ 
and $S_3'$ and on their negative normal bundles. According to Lemma 
\ref{L:action}, the action of $r^*$ on both $H^0(S_2')$ and $H^0(S_3')$ 
is trivial.

\begin{lemma}\label{L:thom1}
The involution $r$ acts as minus identity on the normal bundle of $S_2'
\subset \M'$ and, in particular, on its negative normal bundle.
\end{lemma}

\begin{proof}
This will follow from the fact that $r: \M \to \M$ given by the formula
(\ref{E:res}) is an involution and that its fixed point 
set coincides with $S_2$. To show the latter, suppose that $[A_1,B_1,A_2,
B_2] = [-A_1,B_1,A_2,B_2]$. Then there is $u \in SU(2)$ such that $u A_1 
= - A_1 u$ and $u$ commutes with $B_1$, $A_2$ and $B_2$. Since $(A_1,B_1,
A_2,B_2)$ is an irreducible representation, we conclude that $u^2 = 1$ or 
$u^2 = -1$. The former cannot happen because then $-A_1 = A_1$ and $A_1 = 
0$, a contradiction with $A_1 \in SU(2)$. Therefore, we have $u^2 = - 1$, 
and $u = \pm i$ after conjugation. The fact that $u = \pm i$ commutes 
with $B_1$, $A_2$ and $B_2$ means that all three of them are unit complex 
numbers. In fact, one can conjugate by $j$ if necessary to make $B_1 = 
e^{i\beta}$ with $\beta \in [0,\pi]$, perhaps at the expense of changing 
the sign of $u$. Conjugate further by a unit complex number to make
$A_1 = a + bj$ with real non-negative $b$, while preserving $u$, $B_1$, 
$A_2$ and $B_2$. The relation $[A_1,B_1] = [A_1,B_1]\cdot [A_2,B_2] = -1$ 
then tells us that, up to conjugation, $A_1 = j$ and $B_1 = i$. Therefore,
the only fixed points of $r: \M \to \M$ are of the form $[j,i,z,w]$, with 
$z$ and $w$ some unit complex numbers. This is exactly $S_2 \subset \M$.
\end{proof}

Since the rank of the negative normal bundle to $S_2'$ is one, $r$ 
reverses orientation on the fiber hence $r^*$ anticommutes with the first 
arrow in (\ref{E:thom1}). On the other hand, $S'_3$ is the absolute 
maximum hence its negative normal bundle is the same as its normal bundle. 
According to Lemma \ref{L:action}, the action of $r$ on $S'_3$ is 
orientation preserving. Since $r: \M' \to \M'$ preserves orientation 
(which follows for example from Lemma \ref{L:thom1}) it must be 
orientation preserving on the negative normal bundle of $S_3'$. Therefore, 
$r^*$ commutes with the last arrow in (\ref{E:thom1}). In conclusion, 
$d_1: E_1^{1,0} \to E_1^{2,0} $ anticommutes with $r^*$ and thus must 
vanish.

Finally, let us consider the differential $d_1: E_1^{0,1}\to E_1^{1,1}$.
It is represented as the composition 
\begin{equation}\label{E:thom2}
H^1 (S_1')\; \to\; H^1 (U_1',U_0')\; \to\; H^1 (U_1')\;\to\; 
H^2 (U_2',U_1')\; \to\; H^1 (S_2').
\end{equation}
As before, the two arrows in the middle commute with $r^*$, and the first 
and the last arrows are Thom isomorphisms. The action of $r$ on $S_1'$ is 
by the $180^{\circ}$ rotation hence the induced action on $H^1 (S_1')$ is 
trivial. Since the negative normal bundle of $S_1'$ is zero dimensional, 
there is no Thom class to worry about and we readily conclude that $r^*$ 
commutes with the first arrow in (\ref{E:thom2}). Concerning the last 
arrow, we already know that $r^*$ anticommutes with the Thom class of the 
negative normal bundle of $S_2'$. On the other hand, $r$ acts as the 
identity on $S'_2$ and hence as the identity on $H^1(S'_2)$. This implies 
that $r^*$ anticommutes with the last arrow and hence with the composition 
(\ref{E:thom2}). Thus $d_1: E_1^{0,1} \to E_1^{1,1}$ vanishes. 


\section{Differential $d_2$}
The results of the previous section imply that $E_2 = E_1$ hence all that 
remains to do is compute the differential $d_2: E_2^{0,1}\to E_2^{2,0}$.
We will show that the edge homomorphism $i^*: H^1(\M') \to H^1 (S_1')$ in 
the spectral sequence induced by the inclusion $i: S_1'\to \M'$ is 
surjective. This will imply that $d_2 = 0$ because a generator of 
$E_2^{0,1} = H^1(S'_1) = \Z$ must survive in the $E_{\infty}$ term of the 
spectral sequence, hence it must be in the kernel of $d_2$. Note that a
similar argument could be used to show vanishing of $d_1: E_1^{0,1} \to 
E_1^{1,1}$ above.

Remember that $S'_1$ consists of the quadruples $[e^{i\phi},-1, i, j] \in 
\M'$. Consider the subset $\R' \subset \M'$  which consists of quadruples 
$[1,B_1,A_2,B_2]$ fixed by the involution $\sigma^*: \M' \to \M'$. 

\begin{lemma}
$\R' \subset \M'$ is a smoothly embedded 2--sphere which intersects $S_1' 
\subset \M'$ transversely in exactly one point.
\end{lemma}

\begin{proof}
The relation $[1,B_1] \cdot [A_2,B_2] = [A_2,B_2] = -1$ on the points of
$\R'$ implies that, up to conjugation, $A_2 = i$ and $B_2 = j$. The fact
that $[1,B_1,i,j] \in \R'$ is fixed by $\sigma^*$ means that $\sigma^*
([1,B_1,i,j]) = [1, B_1^{-1},i,-j] = [1, i\,B_1^{-1} i^{-1},i,j]$ (we 
used conjugation by $i$ in the last equality). Therefore, $\R'$ is 
parametrized by $B_1 \in SU(2)$ such that $B_1\,i = i\,B_1^{-1}$.  These
are precisely the unit quaternions with no $i$ component; they obviously 
comprise an embedded 2--sphere in $\M'$. The intersection $\R'\,\cap\,S_1'$
consists of just one point, $[1,-1,i,j]$, and it is obviously transversal.
\end{proof}

Now, the lemma implies that the natural generator of $H^1 (S'_1) = \Z$ 
is the image under $i^*$ of the Poincar\'e dual of $\R'$. This shows 
that $i^*: H^1(\M') \to H^1 (S'_1)$ is surjective and thus completes the 
argument.


\section{Proof of Theorem \ref{T:main}}\label{S:proof}
The Morse--Bott spectral sequence associated with the function  $f':  \M' 
\to \mathbb R$ converges to $H^* (\M') = H^* (\N')$.  As we showed in the 
last two sections, all the differentials in this spectral sequence vanish 
and therefore it collapses at the $E_1$ term, $E_1 = E_{\infty}$. 
This completes the proof in the case when the puncture belongs to the left 
circle of $F'$ as shown in Figure \ref{fig2}.

If the puncture belongs to the right circle of $F'$, the involution 
$\sigma^*: \M \to \M$ is given by the formula
\[
\sigma^* [A_1,B_1,A_2,B_2] = [-B_1\,A_1\,B_1^{-1},\;B_1^{-1},\;
B_2\,A_2\,B_2^{-1},\;B_2^{-1}],
\]
and Theorem \ref{T:main} follows by simply renaming the variables. In 
the remaining case, when $p \in F'$ belongs to the middle circle, the 
involution $\sigma^*: \M \to \M$ is given by 
\[
\sigma^* [A_1,B_1,A_2,B_2] = [B_1\,A_1\,B_1^{-1},\;B_1^{-1},\;
B_2\,A_2\,B_2^{-1},\;B_2^{-1}].
\]
The above proof goes through with minimal changes, which can be safely 
left to the reader.

\begin{remark}\label{R:rings}
The regular neighborhood of $S_1'\,\cup\,\R' = S^1\,\vee\,S^2$ in $\M'$ 
is a punctured $S^1 \times S^2$.  This implies that $\M'$ splits into a 
connected sum, one of the factors being $S^1 \times S^2$. In particular, 
the isomorphism $H^* (\N') = H^* (S^1 \times S^1 \times S^1)$ of Theorem 
\ref{T:main} is not a ring isomorphism.
\end{remark}



\begin{thebibliography}{10}

\bibitem{AB} 
M.~Atiyah and R.~Bott, {\em The Yang-Mills equations over Riemann 
surfaces}, Phil. Tran. R. Soc. A \textbf{308} (1982), 523--615

\bibitem{Da} 
G.~Daskalopoulos, {\em The topology of the space of stable bundles on a 
compact Riemann surface},  J. Differential Geom. \textbf{36} (1992),  
699--746

\bibitem{D} 
J.~Duistermaat, {\em Convexity and tightness for restrictions of 
Hamiltonian functions to fixed point sets of an antisymplectic 
involution}, Trans. Amer. Math. Soc. \textbf{275} (1983), 417--429

\bibitem{F} 
T.~Frankel, {\em Fixed points on K\"ahler manifolds}, Ann. of Math. 
\textbf{70} (1959), 1--8

\bibitem{G} 
W.~Goldman, {\em The symplectic nature of fundamental groups of 
surfaces}, Adv. Math. \textbf{54} (1984), 200--225

\bibitem{HL} 
N.-K.~Ho and C.-C.~Liu, {\em Yang-Mills connections on orientable and 
nonorientable surfaces}. Preprint {\tt arXiv:0707.0258}

\bibitem{NS}
M.~Narasimhan, C.~Seshadri,
{\em Stable and unitary vector bundles on a compact Riemann surface},
Ann. of Math. \textbf{82} (1965), 540--567

\bibitem{T}
M.~Thaddeus, {\em A perfect Morse function on the moduli space of flat 
connections}, Topology \textbf{39} (2000), 773--787

\bibitem{W1} 
S.~Wang, {\em Classification of real moduli spaces over genus $2$ 
curves},  Geom. Dedicata  {\bf57} (1995), 207--215

\bibitem{W2} 
S.~Wang, {\em Orientability of real parts and spin structures}, 
JP J. Geom. Topol. {\bf 7} (2007), 159--174

\end{thebibliography}
\end{document}